\documentclass[a4paper,12pt]{amsart}

\usepackage{t1enc}
\def\R{\mathbb{R}}

\def\al{\alpha}

\def\de{\delta}

\def\om{\omega}

\def\Up{\Upsilon}

\newcommand{\newc}{\newcommand}

\newtheorem{theorem}{Theorem}[section]

\newtheorem{corollary}[theorem]{Corollary}
\theoremstyle{remark}

\newtheorem*{remark*}{\rm\bf Remark}

\usepackage{amssymb,stmaryrd}
\usepackage{amscd}

\newcommand{\nd}{\nabla}

\def\P{\mathrm{P}}

\newc{\aR}{\mbox{\boldmath{$ R$}}}
\newc{\aS}{\mbox{\boldmath{$ S$}}}
\newc{\aDeR}{\mbox{\boldmath{$ U$}}_B{}^P{}_C{}^Q}
\newc{\aDe}{\mbox{\boldmath$ \Delta$}}
\newc{\aNd}{\mbox{\boldmath$ \nabla$}}

\newc{\aK}{\mbox{\boldmath{$ K$}}}
\newc{\aL}{\mbox{\boldmath{$ L$}}}

\def\sideremark#1{\ifvmode\leavevmode\fi\vadjust{\vbox to0pt{\vss% the remark
 \hbox to 0pt{\hskip\hsize\hskip1em%                          will appear only
 \vbox{\hsize3cm\tiny\raggedright\pretolerance10000%          on the side
 \noindent #1\hfill}\hss}\vbox to8pt{\vfil}\vss}}}%
\author{Matthew Randall}
\address{Mathematical Sciences Institute\\ 
Australian National University\\ 
Canberra, ACT 0200\\
 Australia}
\email{matthew.randall@anu.edu.au}
\title{Local obstructions to projective surfaces admitting skew-symmetric Ricci tensor}

\begin{document}


\begin{abstract}
The equation determining whether a projective structure admits a connection in its given projective class that has skew-symmetric Ricci tensor is an overdetermined system of semi-linear partial differential equations which we call the projective Einstein-Weyl (pEW) equation. In 2-dimensions, we give local obstructions for projective surfaces to admit such a connection in its projective class. The obstructions are the resultants of polynomial equations that have to be satisfied for there to admit any pEW solution. 
\end{abstract}

\maketitle

\pagestyle{myheadings}
\markboth{Randall}{Local obstructions to projective surfaces admitting skew-symmetric Ricci tensor}

\section{Introduction}
A projective structure on a smooth manifold is an equivalent class of torsion-free affine connections that have the same unparameterised curves as geodesics. For projective manifolds $(M,[\nd])$, a natural geometric problem is to find an affine connection in the projective class with Ricci tensor identically zero. This can be seen as a projective analogue of Einstein's equation and can be reformulated as solving an overdetermined system of linear partial differential equations, known as the projective Ricci-flat equation in \cite{projectiveBGG}. A generalisation of the Ricci-flat condition is to find an affine connection in the projective class with the symmetric part of the Ricci tensor identically zero, i.e.\ the Ricci tensor is skew-symmetric. The overdetermined system of partial differential equations associated to this condition, which we call the projective Einstein-Weyl (pEW) equation, now becomes semi-linear. In this paper we derive local obstructions to existence of solutions to the pEW equation on 2-dimensional projective manifolds. 
The obstructions are the resultants of polynomials with coefficients given by invariants of the projective structure. Computing these obstructions allows us to distinguish projective surfaces that admit skew-symmetric Ricci tensor from those that do not. In dimension $2$, affine structures with skew-symmetric Ricci-tensor are of interest, as investigated in \cite{Derd} and \cite{AffOss}. We first set up the closed system for the pEW equation in $2$-dimensions, then derive the algebraic constraints that give rise to the obstructions. We conclude with two examples for which the obstructions do not vanish. Abstract index notation as explained in \cite{SStime} is used throughout the paper to describe tensors on the manifold.

\section{Projective differential geometry and the pEW equation}
In this section we review the projective differential geometry needed for the results. More details can be found in \cite{east}.
Let $(M^n,[\nd])$ be a smooth manifold with a smooth projective structure with $n \geq 2$. Two torsion-free affine connections are projectively equivalent if and only if they have the same geodesics up to reparametrization. Equivalently, 
\[
\hat \nd_a\om_b=\nd_a\om_b-\Up_a\om_b-\Up_b\om_a
\]
for some $1$-form $\Up_a$.
The curvature of any affine connection $\nd \in [\nd]$ decomposes as follows:
\begin{equation}\label{curvdecomp}
R_{ab}{}^c{}_d=W_{ab}{}^c{}_d+\de_{a}{}^c\P_{bd}-\de_b{}^c\P_{ad}-2\P_{[ab]}\de^c{}_d
\end{equation}
where $R_{ab}=R_{ca}{}^c{}_b$ is the Ricci curvature and $\P_{ab}=\frac{n}{n^2-1}R_{ab}+\frac{1}{n^2-1}R_{ba}$ is the projective rho tensor. The Ricci tensor (and projective rho tensor $\P_{ab}$) is not necessarily symmetric. The totally trace-free part of $R_{ab}{}^c{}_d$ denoted by $W_{ab}{}^c{}_d$ is called the projective Weyl tensor and is an invariant of the projective structure. In dimension $2$, the Weyl tensor $W_{ab}{}^c{}_d$ vanishes by symmetry considerations. The skew-symmetric part of the projective rho tensor $\P_{[ab]}$ is called the Faraday 2-form $F_{ab}$. It is always closed. The condition of a projective structure admitting skew-symmetric Ricci tensor is equivalent to finding an affine connection $\nd \in [\nd]$ such that the symmetric part of its projective rho tensor can be made to vanish, i.e.\ $\P_{(ab)}=0$.
In this case $\P_{ab}=\P_{[ab]}=F_{ab}$ and the formula (\ref{curvdecomp}) reduces to  
\[
R_{ab}{}^c{}_d=W_{ab}{}^c{}_d+\de_{a}{}^cF_{bd}-\de_b{}^cF_{ad}-2F_{ab}\de^c{}_d.
\]
Under a projective transformation, the rho tensor rescales as 
\[
\hat \P_{ab}=\P_{ab}-\nd_a\Up_b+\Up_a\Up_b. 
\]
The above condition is then equivalent to solving for $0=\hat \P_{(ab)}=\P_{(ab)}-\nd_{(a}\Upsilon_{b)}+\Upsilon_a\Upsilon_b$. Making the substitution $\al_a=-\Up_a$ gives the pEW equation
\begin{equation}\label{pEW}
\nd_{(a}\al_{b)}+\al_a\al_b+\P_{(ab)}=0.
\end{equation}
The $1$-form $\al_a$ changes by a gauge $\al_a \mapsto \al_a+\Up_a$ under projective transformations. The pEW equation is a projectively invariant overdetermined system of semi-linear PDEs, and specialises to the projective Ricci-flat (or projective to Einstein) equation when $F_{ab}=0$ (i.e.\ $\al_a$ is locally exact). In $2$-dimensions, we shall fix a volume form $\epsilon^{ab}$ with $\epsilon^{ab}\epsilon_{ab}=2!=2$, and we set our convention so that $\epsilon^{ac}\epsilon_{bc}=\delta^a_b$. The preferred choice of affine connection $\nd \in [\nd]$ that we work with is one such that $\nd_a\epsilon_{bc}=0$, and there is no obstruction even globally in choosing that. As a consequence of having a parallel volume form, the projective rho tensor $\P_{ab}=\P_{(ab)}$ is now symmetric. Also in dimension $2$, the projective Cotton-York tensor defined by $Y_{abc}:=\nd_a\P_{bc}-\nd_b\P_{ac}$ is projectively invariant. Using the volume form to dualise, we can write $Y_{abc}=\frac{1}{2}\epsilon_{ab}Y_c$, where $Y_a=\epsilon^{bc}Y_{bca}$. The vanishing of $Y_a$ characterises projectively flat surfaces.

\subsection{Deriving the closed system for pEW on projective surfaces}
A common procedure to treat equations such as (\ref{pEW}) is through prolongation \cite{prolong}. This involves expressing first derivatives of the dependent variables in terms of the variables themselves. Introduce $F_{ab}=\nd_{[a}\al_{b]}$ as the extra dependent variable. Using the volume form to dualise, we can write $F_{ab}=\frac{1}{2}\epsilon_{ab}F$ where $F=\epsilon^{ab}F_{ab}$ is a projective scalar density of projective weight $-3$ (see \cite{east} for the definition of projective weight). We can rewrite (\ref{pEW}) as
\begin{equation}\label{pEW2}
\nd_a\al_b+\al_a\al_b+\P_{ab}=\frac{1}{2}\epsilon_{ab}F.
\end{equation}
Differentiating (\ref{pEW2}), we find that
\begin{align}\label{2pew01}
\nd_aF=&-3F\al_a-Y_a
\end{align}
is a consequence of the original equation. Equations (\ref{pEW2}) and (\ref{2pew01}) form the first order closed system for the pEW equation, and from this we can derive algebraic constraints for (\ref{pEW}) to hold by further differentiating the system. 

\section{Statement of results}
In the flat case when $Y_a=0$, we necessarily have $F=0$ by (\ref{2pew02}). The $1$-form $\al_a$ is therefore exact, and equation (\ref{pEW}) specialises to the projective Ricci-flat equation. We shall now restrict our attention to non-flat projective surfaces, that is one with $Y_a$ non-zero. This ensures that $F\neq 0$. 
Let $\text{Res}(P(t),Q(t))$ be the resultant of polynomials $P(t)$ and $Q(t)$ in the single variable $t$. $\text{Res}(P(t),Q(t))=0$ is a necessary and sufficient condition for $P(t)$ and $Q(t)$ to share a common root. It turns out in deriving the constraint equations for (\ref{pEW}) to hold we have to distinguish between the cases where a certain projective invariant $\rho:=Y_aW^a$ vanishes or not. The quantity $W^a$ will be defined later in section \ref{proof1rhonot0}.   
\begin{theorem}\label{1rhonot0}
Let $(M^2,[\nd])$ be a projective surface with $\rho \neq 0$. Suppose $M^2$ admits a solution to (\ref{pEW}). Then there exists polynomials $P_1(t)$, $P_2(t)$, $P_3(t)$ in the single variable $t$ with coefficients given by the invariants of the projective structure such that when $t=F$, 
\[
P_1(F)=P_2(F)=P_3(F)=0
\]
must hold.   
\end{theorem}
The polynomial constraints $P_1(F)=P_2(F)=P_3(F)=0$ are explicitly computed in section \ref{proof1rhonot0}. As a corollary, we obtain local obstructions for there to be solutions for (\ref{pEW}).
\begin{corollary}\label{rhonot0}
Let $(M^2,[\nd])$ be a projective surface with $\rho \neq 0$. Suppose $M^2$ admits a solution to (\ref{pEW}). Then 
\[
\text{ Res}(P_1(F), P_2(F))=\text{ Res}(P_2(F), P_3(F))=\text{ Res}(P_1(F), P_3(F))=0
\]  
must hold.
\end{corollary}
\begin{theorem}\label{rho0}
Let $(M^2,[\nd])$ be a projective structure with $\rho=0$. Suppose $M^2$ admits a solution to (\ref{pEW}). Then 
\begin{equation}\label{obst2}
kY^a\nd_am-m(Y^a\nd_ak+k\phi-6m)=0
\end{equation}
must hold, where $k, m, \phi$ are quantities obtained from the projective structure on $M^2$ to be defined later in sections \ref{proof1rhonot0} and \ref{proofrho0}. 
\end{theorem}

\section{Proof of theorem \ref{1rhonot0}}\label{proof1rhonot0}
We shall now derive the polynomial equations $P_1(F)=0$, $P_2(F)=0$, $P_3(F)=0$ that arise for (\ref{pEW}) to hold. The polynomials $P_1(t)$, $P_2(t)$, $P_3(t)$ are obtained by replacing $F$ with the indeterminate $t$. Using the volume form to raise and lower indices, we have $Y^a=\epsilon^{ab}Y_b$, $Y_b=Y^a\epsilon_{ab}$.
Differentiating (\ref{2pew01}) and using the closed system for the pEW equation gives 
\begin{align}\label{pEWFF}
\nd_a\nd_bF=&9F\al_a\al_b+3Y_a\al_b-\frac{3}{2}\epsilon_{ab}F^2+3\al_a\al_bF+3\P_{ab}F-\nd_aY_{b}
\end{align}
and skewing with the volume form $\epsilon^{ab}$ gives
\begin{equation}\label{2pew02}
3\al_aY^a+3F^2+\nd_aY^a=0 
\end{equation}
as the first constraint of the system. 
Let
\[
\phi=2\nd_aY^a
\]
and we can rewrite equation (\ref{2pew02}) as
\begin{equation}\label{2pew05}
\al_aY^a=-F^2-\frac{\phi}{6}.
\end{equation}
Introduce $W^a=Y^b\nd_bY^a-\frac{2\phi}{3}Y^a$ a projectively invariant vector (of projective weight $-12$). 
Differentiating (\ref{2pew05}) and using (\ref{pEW2}), (\ref{2pew01}) and (\ref{2pew05}), we find that
\[
\al_aW^a=-5F^4+\frac{5}{36}\phi^2+\P_{ab}Y^aY^b-\frac{Y^a\nd_a\phi}{6}.
\]
Let
\begin{equation*}
\ell=\frac{5\phi^2}{12}+3\P_{ac}Y^aY^c-\frac{Y^a\nd_a\phi}{2}.
\end{equation*}
We thus have
\begin{equation}\label{2pew06}
\al_aW^a=-5F^4+\frac{\ell}{3}
\end{equation}
as our second constraint.  
We can now solve for $\al_a$ assuming $\rho=Y_aW^a \neq 0$. It is given by
\begin{align}\label{newal}
\al_a=&\frac{1}{3\rho}\left(\frac{\phi}{2}+3F^2\right)W_a-\frac{1}{3\rho}\left(15F^4-\ell\right)Y_a.
\end{align}
Substituting (\ref{newal}) back into equation (\ref{pEW}) yields further constraints on $F$. The first polynomial constraint $P_1(F)=0$ comes from computing $F=\nd_a\al^a$ using (\ref{newal}). It is given by
\begin{align*}
P_1(F)=&-90F^6+15\left(\frac{Y^a\nd_a\rho}{\rho}-\frac{5\phi}{2}\right)F^4-\left(\frac{3W^a\nd_a\rho}{\rho}+6\ell-3\nd_aW^a\right)F^2\\
&-9\rho F+\left(\frac{W^a\nd_a\phi}{2}+\frac{\phi}{2}\nd_aW^a+Y^a\nd_a\ell+\frac{\phi \ell}{2}-\frac{\phi W^a\nd_a\rho}{2\rho}-\ell\frac{Y^a\nd_a\rho}{\rho}\right)\\
=&0.
\end{align*}   
It can be verified that the coefficients appearing in the polynomial $P_1(t)$ are all projectively invariant. For example, under projective rescalings, the coefficient of the degree $4$ term in $P_1(t)$ transforms as follows:
\begin{align*}
15\widehat{\left(\frac{Y^a\nd_a\rho}{\rho}-\frac{5\phi}{2}\right)}=15\left(\frac{Y^a\hat\nd_a\rho}{\rho}-\frac{5\hat\phi}{2}\right)=&15\left(\frac{Y^a\nd_a\rho-15\Up_a Y^a \rho}{\rho}-\frac{5(\phi-6\Up_aY^a)}{2}\right)\\
=&15\left(\frac{Y^a\nd_a\rho}{\rho}-15\Up_a Y^a-\frac{5\phi}{2}+15\Up_aY^a\right)\\
=&15\left(\frac{Y^a\nd_a\rho}{\rho}-\frac{5\phi}{2}\right).
\end{align*}
It is therefore projectively invariant. The second and third polynomial constraints come from substituting (\ref{newal}) back into (\ref{pEW}) and contracting with $W^aW^b$ and $W^aY^b$ respectively.
Another possible contraction with $Y^aY^b$ yields an equation that is identically zero.  
Evaluating $W^aW^b\nd_a\al_b+\al_a\al_bW^aW^b+\P_{ab}W^aW^b=0$ for $\al_a$ in (\ref{newal}) gives
\begin{align*}
P_2(F)=&-275F^8+\left( \frac{-5W^eY^d\nd_eW_d}{\rho}+\frac{50 \ell}{3} \right)F^4+20\rho F^3+\left(\frac{W^eW^a\nd_eW_a}{\rho}\right)F^2\\
&+\frac{\phi W^eW^a\nd_eW_a}{6\rho}+\P_{ea}W^eW^a+\frac{\ell^2}{9}+\frac{\ell W^eY^d\nd_eW_d}{3\rho}+\frac{W^e\nd_e\ell}{3}\\
=&0,
\end{align*}
while evaluating $W^aY^b\nd_{(a}\al_{b)}+\al_a\al_bW^aY^b+\P_{ab}W^aY^b=0$ gives
\begin{align*}
P_3(F)=&-40F^6+\left(\frac{-5Y^a(W^e\nd_eY_a+Y^e\nd_aW_e)}{2\rho}-\frac{25\phi}{6}\right)F^4\\
&+\left(\frac{2\ell}{3}+\frac{W^e(W^d\nd_eY_d+Y^d\nd_dW_e)}{2\rho}\right)F^2\\
&+\rho F-\frac{W^e\nd_e\phi}{12}+\frac{\phi W^bW^e\nd_eY_b}{12\rho}\\
&+\frac{\ell W^eY^a\nd_eY_a}{6\rho}+\P_{ae}W^eY^a+\frac{Y^e\nd_e\ell}{6}+\frac{\phi W^eY^a\nd_aW_e}{12\rho}-\frac{\ell \phi}{18}+\frac{\ell Y^bY^a\nd_aW_b}{6\rho}\\
=&0. 
\end{align*}
The coefficients of the polynomials $P_2(t)$ and $P_3(t)$ are also all projectively invariant. 
This concludes the proof of Theorem \ref{1rhonot0}. We shall now explain a more concise way of extracting the obstructions. 
\subsection{Concise way of extracting obstructions}\label{binobs}
We can eliminate the single odd degree term so that even degree terms remain in $P_1(t)$, $P_2(t)$, $P_3(t)$. Namely, define
\begin{align*}
Q_1(t^2)=&-20t^2P_3(t)+P_2(t), \qquad Q_2(t^2)=-9P_3(t)-P_1(t),\\
Q_3(t^2)=&-\frac{20}{9}t^2P_1(t)-P_2(t).
\end{align*}
Then the three polynomials $Q_1(t^2)$, $Q_2(t^2)$, $Q_3(t^2)$ will be quartic, cubic and quartic polynomials of $t^2$ respectively since only even powers of $t$ remain. This allows the obstructions to be extracted easily since now the resultant of any of these $2$ polynomials will at most be the determinant of a $8$ by $8$ matrix.  
Supposing that
\begin{align*}
P_1(t)=&-90t^6+a_1t^4+a_2t^2-9\rho t+a_3,\\
P_2(t)=&-275t^8+b_1t^4+20\rho t^3+b_2 t^2+b_3,\\
P_3(t)=&-40t^6+c_1 t^4+c_2 t^2+\rho t+c_3,
\end{align*}
where
\begin{align*}
a_1=&15\left(\frac{Y^a\nd_a\rho}{\rho}-\frac{5\phi}{2}\right)\\
a_2=&-\left(\frac{3W^a\nd_a\rho}{\rho}+6\ell-3\nd_aW^a\right)\\
a_3=&\left(\frac{W^a\nd_a\phi}{2}+\frac{\phi}{2}\nd_aW^a+Y^a\nd_a\ell+\frac{\phi \ell}{2}\right)-\frac{\phi W^a\nd_a\rho}{2\rho}-\ell\frac{Y^a\nd_a\rho}{\rho}\\
b_1=&\frac{-5W^eY^d\nd_eW_d}{\rho}+\frac{50 \ell}{3} \\
b_2=&\frac{W^eW^a\nd_eW_a}{\rho}\\
b_3=&\frac{\phi W^eW^a\nd_eW_a}{6\rho}+\P_{ea}W^eW^a+\frac{\ell^2}{9}+\frac{\ell W^eY^d\nd_eW_d}{3\rho}+\frac{W^e\nd_e\ell}{3}\\
c_1=&\frac{-5Y^a(W^e\nd_eY_a+Y^e\nd_aW_e)}{2\rho}-\frac{25\phi}{6}\\
c_2=&\frac{2\ell}{3}+\frac{W^e(W^d\nd_eY_d+Y^d\nd_dW_e)}{2\rho}\\
c_3=&-\frac{W^e\nd_e\phi}{12}+\frac{\phi W^bW^e\nd_eY_b}{12\rho}+\frac{\ell W^eY^a\nd_eY_a}{6\rho}+\P_{ae}W^eY^a+\frac{Y^e\nd_e\ell}{6}+\frac{\phi W^eY^a\nd_aW_e}{12\rho}\\
&-\frac{\ell \phi}{18}+\frac{\ell Y^bY^a\nd_aW_b}{6\rho},\\
\end{align*}
then a computation gives
\begin{align*}
Q_1(X=t^2)=&525 X^4-20c_1 X^3+(b_1-20c_2)X^2+(b_2-20c_3)X+b_3\\
Q_2(X=t^2)=&450 X^3-(9c_1+a_1)X^2-(9c_2+a_2)X-(9c_3+a_3)\\
Q_3(X=t^2)=&475 X^4-\frac{20}{9}a_1 X^3-\left(\frac{20}{9}a_2+b_1\right)X^2-\left(\frac{20}{9}a_3+b_2\right)X-b_3.
\end{align*}
The local obstructions are therefore
\begin{align*}
&Q_{12}=\\
&\begin{vmatrix}
525 & -20c_1 & (b_1-20c_2)& (b_2-20c_3)& b_3&0&0 \\
0& 525 & -20c_1 & (b_1-20c_2)& (b_2-20c_3)& b_3 &0 \\
0 &0 &525 & -20c_1 & (b_1-20c_2)& (b_2-20c_3)& b_3 \\
450 &-(9c_1+a_1)&-(9c_2+a_2)&-(9c_3+a_3)&0&0&0\\
0&450 &-(9c_1+a_1)&-(9c_2+a_2)&-(9c_3+a_3)&0&0\\
0&0&450 &-(9c_1+a_1)&-(9c_2+a_2)&-(9c_3+a_3)&0\\
0&0&0&450 &-(9c_1+a_1)&-(9c_2+a_2)&-(9c_3+a_3)\\
\end{vmatrix},
\end{align*}
\begin{align*}
&Q_{23}=\\
&\begin{vmatrix}
450 &-(9c_1+a_1)&-(9c_2+a_2)&-(9c_3+a_3)&0&0&0\\
0 &450 &-(9c_1+a_1)&-(9c_2+a_2)&-(9c_3+a_3)&0&0\\
0 & 0&450 &-(9c_1+a_1)&-(9c_2+a_2)&-(9c_3+a_3)&0\\
0 & 0&0&450 &-(9c_1+a_1)&-(9c_2+a_2)&-(9c_3+a_3)\\
475 &-\frac{20}{9}a_1 &-\left(\frac{20}{9}a_2+b_1\right)&-\left(\frac{20}{9}a_3+b_2\right)&-b_3&0&0\\
0&475 &-\frac{20}{9}a_1 &-\left(\frac{20}{9}a_2+b_1\right)&-\left(\frac{20}{9}a_3+b_2\right)&-b_3&0\\
0&0&475 &-\frac{20}{9}a_1 &-\left(\frac{20}{9}a_2+b_1\right)&-\left(\frac{20}{9}a_3+b_2\right)&-b_3\\
\end{vmatrix},
\end{align*}and
\begin{align*}
&Q_{13}=\\
&\begin{vmatrix}
525 & -20c_1 & (b_1-20c_2)& (b_2-20c_3)& b_3&0&0&0 \\
0& 525 & -20c_1 & (b_1-20c_2)& (b_2-20c_3)& b_3 &0&0 \\
0 &0 &525 & -20c_1 & (b_1-20c_2)& (b_2-20c_3)& b_3&0 \\
0 &0 & 0& 525 & -20c_1 & (b_1-20c_2)& (b_2-20c_3)& b_3\\
475 &-\frac{20}{9}a_1 &-\left(\frac{20}{9}a_2+b_1\right)&-\left(\frac{20}{9}a_3+b_2\right)&-b_3&0&0&0\\
0&475 &-\frac{20}{9}a_1 &-\left(\frac{20}{9}a_2+b_1\right)&-\left(\frac{20}{9}a_3+b_2\right)&-b_3&0&0\\
0&0&475 &-\frac{20}{9}a_1 &-\left(\frac{20}{9}a_2+b_1\right)&-\left(\frac{20}{9}a_3+b_2\right)&-b_3&0\\
0&0&0&475 &-\frac{20}{9}a_1 &-\left(\frac{20}{9}a_2+b_1\right)&-\left(\frac{20}{9}a_3+b_2\right)&-b_3\\
\end{vmatrix},
\end{align*}
and the knowledge of the values of $a_1, a_2, a_3, b_1, b_2, b_3, c_1, c_2, c_3$ at a point will let us determine the values of the obstruction at that point.

\section{Proof of theorem \ref{rho0}}\label{proofrho0}
In the case where $\rho=0$, the vectors $W^a$ and $Y^a$ are linearly dependent everywhere on $M^2$ and we can express $W^a=fY^a$ (for some projectively invariant density $f$ of weight $-6$). Constraint (\ref{2pew06}) becomes 
\begin{align*}
f(\al_aY^a)=-5F^5+\frac{\ell}{3},
\end{align*}
upon which multiplying (\ref{2pew05}) by $f$ and eliminating $\al_aY^a$ gives the quartic equation  
\begin{equation}\label{cons}
15F^4-3fF^2-\left(\ell+\frac{f \phi}{2}\right)=0.
\end{equation}
Let $h=\ell+\frac{f \phi}{2}$.
Differentiating (\ref{cons}) gives
\begin{equation*}
60F^3(\nd_aF)-6fF(\nd_aF)-3(\nd_af)F^2-\nd_ah=0, 
\end{equation*}
which we contract with $Y^a$ and use (\ref{2pew01}) and (\ref{2pew05}) to get
\begin{align}\label{cons0}
180F^6+(30\phi-18f)F^4-3(\phi f+Y^a\nd_af)F^2-Y^a\nd_ah=0.
\end{align}
Multiplying (\ref{cons}) by $12F^2$ gives
\[
180F^6-36fF^4-12hF^2=0,
\]
which can be used to eliminate the term of degree $6$ in (\ref{cons0}) to obtain
\begin{equation}\label{cons1}
(30\phi+18f)F^4-3(\phi f+Y^a\nd_af-4h)F^2-Y^a\nd_ah=0. 
\end{equation}
Using (\ref{cons}) once more, we can further eliminate the term of degree $4$ in (\ref{cons1}) to obtain the quadratic equation
\begin{equation}\label{quad1}
\left(3\phi f-3Y^a\nd_af+12h+\frac{18}{5}f^2\right)F^2+\frac{6hf}{5}-Y^a\nd_ah+2\phi h=0.
\end{equation}
Let us call
\begin{align*}
k=&3\phi f-3Y^a\nd_af+12h+\frac{18}{5}f^2,\\
m=&\frac{6hf}{5}-Y^a\nd_ah+2\phi h,
\end{align*}
so that (\ref{quad1}) becomes
\begin{equation}\label{cons2}
kF^2+m=0.
\end{equation}
Differentiating (\ref{cons2}) gives 
\[
(\nd_ak)F^2+2Fk\nd_aF+\nd_am=0, 
\]
which contracted into $Y^a$ and again using (\ref{2pew01}) and (\ref{2pew05}) gives
\begin{align*}
6kF^4+(k\phi+Y^a\nd_ak)F^2+Y^a\nd_am=0.
\end{align*}
Now substituting 
\begin{align*}
6kF^4&=-6mF^2,
\end{align*}
which is a consequence of (\ref{cons2}), we obtain
\begin{equation}
(k\phi+Y^a\nd_ak)F^2+Y^a\nd_am-6mF^2=0,
\end{equation}
which we multiply by $k$ to get
\[
(k\phi+Y^a\nd_ak)kF^2+kY^a\nd_am-6mkF^2=0. 
\]
Again we make use of (\ref{cons2}) to get
\begin{align*}
-m(k\phi+Y^a\nd_ak)+kY^a\nd_am+6m^2=kY^a\nd_am-m(Y^a\nd_ak+k\phi-6m)=0,
\end{align*}
which is the desired obstruction (\ref{obst2}).

\section{Examples}
In this section we give 2 different projective structures on $\R^2$ that yield non-vanishing obstructions, one with the projective invariant $\rho \neq 0$ and the other with $\rho=0$.

\subsection{Example with $\rho \neq 0$ and non-vanishing obstruction}
This projective structure on $\R^2$ has the the connection coefficients given by
\[
\Pi^1_{22}=xy,\hspace{12pt} \Pi^2_{11}=-y, \hspace{12pt} \Pi^1_{11}=\Pi^1_{21}=\Pi^2_{12}=\Pi^2_{22}=0.
\]   
We compute the polynomials $P_1(t)$, $P_2(t)$ and $P_3(t)$ at an arbitrary point $p \in \R^2$ where $\rho(p)\neq 0$. Taking $p$ to be given in local coordinates by $(x,y)=(1,1)$, we find $\rho(p)=328$, and the polynomials at the point $p$ are given by
\begin{align*}
P_1(t)&=-90t^6+\frac{185760}{328}t^4-\frac{528608}{328}t^2-2952t-\frac{134912}{328},\\
P_2(t)&=-275t^8+\frac{13774080}{2952}t^4+6560t^3-\frac{601856}{8856}t^2+\frac{523957248}{26568}\\
P_3(t)&=-40t^6+\frac{30960}{328}t^4+\frac{125360}{984}t^2+328t+\frac{31603200}{26568}. 
\end{align*}
Using the software MAPLE, we find that 
\begin{align*}
Q_{12}(p)=&-\frac{1457890459574161592339200000}{1681},\\
Q_{13}(p)=&-\frac{188610437798501965389961756672000000000}{452190681}, \\
Q_{23}(p)=&\frac{1457890459574161592339200000}{1681}.
\end{align*}
These three polynomials therefore do not share a common root, hence the obstructions do not vanish on any open set containing $p$. Since $p$ is chosen arbitrarily, we conclude that this projective structure does not admit any solution to pEW locally.

\subsection{Example with $\rho=0$ and non-vanishing obstruction}
This projective structure on $\R^2$ has the the connection coefficients given by
\[
\Pi^1_{11}=-\frac{x^2}{6}, \hspace{12pt} \Pi^1_{22}=-\frac{x^2}{2},\hspace{12pt} \Pi^2_{21}=\frac{x^2}{6}, \hspace{12pt} \Pi^2_{11}=\Pi^1_{21}=\Pi^2_{22}=0.
\]  
A computation shows that $\rho=0$. With the help of MATLAB, we find that the obstruction is
\begin{align*}
&k(Y^a\nd_am)-m(Y^a\nd_ak+k\phi-6m)\\
=&\frac{32}{75}x^{30}+\frac{144}{25}x^{27}-208x^{24}-\frac{117408}{25}x^{21}+\frac{311104}{5}x^{18}+136208x^{15}+71536x^{12}\\
&-996160x^9+978560x^6-\frac{824000}{3}x^3+27200
\end{align*}
and we conclude that this projective structure does not admit any solution to pEW locally.

\section{Acknowledgements}
The author would like to thank Michael Eastwood and Paul Tod for comments and discussions, and Graham Weir for his observation leading to subsection \ref{binobs}.

\end{document}